\documentclass[a4paper,12pt]{amsart}
\usepackage{amsfonts}
\usepackage{amssymb}
\usepackage{ifthen}
\usepackage{graphicx}
\usepackage[usenames]{color}

\nonstopmode \numberwithin{equation}{section}
\newtheorem{thm}{Theorem}[section]

\newtheorem{cor}{Corollary}[section]
\newtheorem{lem}{Lemma}[section]
\newtheorem{prop}{Proposition}[section]

\newtheorem{claim}{Claim}
\newtheorem{conj}{Conjecture}

\theoremstyle{definition}
\newtheorem{defn}{Definition}[section]
\newtheorem{examp}{Example}[section]
\newtheorem{prob}[equation]{Problem}
\newtheorem{ques}{Question}[section]
\newtheorem{rem}{Remark}[section]

\newcounter {own}
\def\theown {\thesection       .\arabic{own}}

\newenvironment{pf}[1][]{%
 \vskip 3mm
 \noindent
 \ifthenelse{\equal{#1}{}}%
  {{\slshape Proof. }}%
  {{\slshape #1.} }%
 }%
{\qed\bigskip}

\newcounter{alphabet}

\newenvironment{Thm}[1][]{\refstepcounter{alphabet}%
\bigskip%
\noindent%
{\bf Theorem \Alph{alphabet}}%
\ifthenelse{\equal{#1}{}}{}{ (#1)}%
{\bf .} \itshape}{\vskip 8pt}

\newcounter{alphabet2}

\newenvironment{Lem}[1][]{\refstepcounter{alphabet}%
\bigskip%
\noindent%
{\bf Lemma \Alph{alphabet}}%
{\bf .} \itshape}{\vskip 8pt}

\newenvironment{Conj}[1][]{\refstepcounter{alphabet}%
\bigskip%
\noindent%
{\bf Conjecture \Alph{alphabet}}%
{\bf .} \itshape}{\vskip 8pt}

\newcommand{\IC}{{\mathbb C}}

\newcommand{\ID}{{\mathbb D}}

\makeatletter
\@namedef{subjclassname@2020}{%
  \textup{2020} Mathematics Subject Classification}
\makeatother

\def\be{\begin{equation}}
\def\ee{\end{equation}}

\newcommand{\bee}{\begin{enumerate}}
\newcommand{\eee}{\end{enumerate}}

\newcommand{\blem}{\begin{lem}}
\newcommand{\elem}{\end{lem}}
\newcommand{\bthm}{\begin{thm}}
\newcommand{\ethm}{\end{thm}}
\newcommand{\bcor}{\begin{cor}}
\newcommand{\ecor}{\end{cor}}
\newcommand{\beg}{\begin{examp}}
\newcommand{\eeg}{\end{examp}}
\newcommand{\begs}{\begin{examples}}
\newcommand{\eegs}{\end{examples}}
\newcommand{\bdefe}{\begin{defn}}
\newcommand{\edefe}{\end{defn}}
\newcommand{\bprob}{\begin{prob}}
\newcommand{\eprob}{\end{prob}}
\newcommand{\bques}{\begin{ques}}
\newcommand{\eques}{\end{ques}}
\newcommand{\bei}{\begin{itemize}}
\newcommand{\eei}{\end{itemize}}

\newcommand{\bca}{\begin{case}}
\newcommand{\eca}{\end{case}}
\newcommand{\bcl}{\begin{claim}}
\newcommand{\ecl}{\end{claim}}

\newcommand{\bcon}{\begin{conj}}
\newcommand{\econ}{\end{conj}}
\newcommand{\bcons}{\begin{conjs}}
\newcommand{\econs}{\end{conjs}}
\newcommand{\bprop}{\begin{prop}}
\newcommand{\eprop}{\end{prop}}
\newcommand{\br}{\begin{rem}}
\newcommand{\er}{\end{rem}}
\newcommand{\brs}{\begin{rems}}
\newcommand{\ers}{\end{rems}}
\newcommand{\bo}{\begin{obser}}
\newcommand{\eo}{\end{obser}}
\newcommand{\bos}{\begin{obsers}}
\newcommand{\eos}{\end{obsers}}
\newcommand{\bpf}{\begin{pf}}
\newcommand{\epf}{\end{pf}}
\newcommand{\ba}{\begin{array}}
\newcommand{\ea}{\end{array}}
\newcommand{\beq}{\begin{eqnarray}}
\newcommand{\beqq}{\begin{eqnarray*}}
\newcommand{\eeq}{\end{eqnarray}}
\newcommand{\eeqq}{\end{eqnarray*}}

\newcommand{\ra}{\rightarrow}

\newcommand{\ds}{\displaystyle}

\begin{document}

\bibliographystyle{amsplain}

\title[Coefficients estimate of $K$-quasiconformal harmonic mappings]
{On the coefficients estimate of $K$-quasiconformal harmonic mappings}
\author{Peijin Li}
\address{Peijin Li,
 Department of Mathematics,
Hunan First Normal University, Changsha, Hunan 410205, People's Republic of China}
\email{wokeyi99@163.com}

\author{Saminathan Ponnusamy$^*$}
\address{S. Ponnusamy, Department of Mathematics,
Indian Institute of Technology Madras, Chennai-600 036, India. }


\email{samy@iitm.ac.in}

\date{\today}


\thanks{$^*$Corresponding author.}

\begin{abstract}
Recently, the Wang et al. \cite{wwrq} proposed a coefficient conjecture for the family
${\mathcal S}_H^0(K)$ of $K$-quasiconformal harmonic mappings $f = h + \overline{g}$ that are sense-preserving and univalent, where
$h(z)=z+\sum_{k=2}^{\infty}a_kz^k$ and $g(z)=\sum_{k=1}^{\infty}b_kz^k$ are analytic in the unit disk $|z|<1$, and the dilatation $\omega =g'/h'$
satisfies the condition $|\omega(z)| \leq k<1$ for $\ID$, with $K=\frac{1+k}{1-k}\geq 1$.  The main aim of this article is provide an affirmative answer in support of  this conjecture by proving this conjecture for every starlike function
(resp. close-to-convex function) from $\mathcal{S}^0_H(K)$. In addition, we verify this conjecture also for
typically real $K$-quasiconformal harmonic mappings. Also,  we establish sharp  coefficients estimate of convex $K$-quasiconformal harmonic mappings.
 By doing so, our work provides a document in support of the main conjecture of  Wang et al..
\end{abstract}


\maketitle 

\section{Introduction and the main results}\label{csw-sec1}


A complex-valued function $f$ defined on  the unit disk $\mathbb{D}:=\{z\in \mathbb{C}:\, |z|< 1\}$ is called harmonic if it satisfies the
Laplace equation $\Delta f =0$, where $\Delta$ denotes the Laplacian operator
$\Delta  = 4 f_{z\overline{z}} = f_{xx} + f_{yy}.
$
If $f$ is normalized by the conditions $f(0)=0$ and $f_{z}(0)=1$, then $f$ has a canonical representation $f = h + \overline{g}$, where,
\be\label{111}
h(z)=z+\sum_{k=2}^{\infty}a_kz^k~\mbox{ and }~ g(z)=\sum_{k=1}^{\infty}b_kz^k,
\ee
and we denote the class of all such functions by $\mathcal{H}$. A locally univalent function $f= h + \overline{g}\in \mathcal{H}$ is called sense-preserving if the Jacobian $J_f(z)=|h'(z)|^2-|g'(z)|^2$ is positive in $\ID$.

\subsection{Conjecture of Clunie and Sheil-Small}
Let ${\mathcal S}_H$ be the class of harmonic functions $f\in\mathcal{H}$ that are univalent and sense-preserving in $\ID$.
We denote by ${\mathcal S}_H^0$,  the subclass of ${\mathcal S}_H$ with the additional condition $f_{\overline{z}}(0)=0$. Several properties of this class together with its various geometric subclasses were investigated first by Clunie and Sheil-Small \cite{CSS} and later by many researchers (cf. \cite{DuR}).
For the class
$\mathcal{S} = \left\{ f=h+\overline{g} \in \mathcal{S}_H :\, g(z) \equiv 0 ~\mbox{ on }
\mathbb{D}\right\},
$
de Branges \cite{de_Branges} has proved that $|a_n| \leq n$ for all $n \geq 2$, and for each $f\in \mathcal{S}$. This settled the Bieberbach conjecture.

A function $f \in \mathcal{S}_H$ is said to be starlike (resp. convex, close-to-convex, convex in the direction $\theta$) in $\ID$ if the range $f(\ID)$ is
starlike with respect to $0$ (resp. convex, close-to-convex, convex in a direction $\theta$), see \cite{CSS,DuR,PonRasi2013, SS}.
The class of all starlike functions (resp. convex, close-to-convex) $f \in \mathcal{S}_H$  is denoted by ${\mathcal S}_H^{*}$
(resp. ${\mathcal C}_H$, ${\mathcal K}_H$).
One can similarly define ${\mathcal S}_H^{*0}$, ${\mathcal C}_H^0$ and ${\mathcal K}_H^0$ in the usual manner (with the assumption that
$b_1=g'(0)=f_{\overline{z}}(0)=0$). A domain $D$ in $\IC$ is called starlike with respect to the $0$ if  $w\in D$ implies $tw\in D$ for all $t\in [0,1]$. Similarly,
a domain $D$  is called convex if  $w,w'\in D$ implies $(1-t)w+tw'\in D$ for all $t\in [0,1]$.

Recall that a domain $D$ in $\mathbb{C}$ is called convex in the direction
$\theta$ $(0\leq \theta < \pi)$ if every line parallel to the line through $0$ and $e^{i\theta}$
has a connected or empty intersection with $D$.
A domain $D\subset\mathbb{C}$ is close-to-convex if its complement $D^c:=\mathbb{C} \setminus D$ of $D$ can be written as a union of non-intersecting rays that go to infinity.


A function $f$ analytic in $\ID$ is said to be {\it typically real} if, $f(z)$ is  real if and only if $z$ is real.
Let ${\mathcal T}$ denote the class of analytic typically real functions $\varphi$ for which $\varphi(0)=0$ and $\varphi'(0)=1$.
The harmonic typically real functions are defined in a similar way. A complex-valued function $f$ harmonic in $\ID$ is said to be {\it typically real} if $f(z)$ is real if and only if $z$ is real. The class ${\mathcal T}_H$ consists of all sense-preserving typically real harmonic functions $f=h+\overline{g}$ with $h(0)=g(0)=0$, $|h'(0)|=1$, and $f(r)>0$ for $0<r<1$. The subclass of ${\mathcal T}_H$ with $g'(0)=0$ is denoted by ${\mathcal T}_H^0$.

In particular, Clunie and Sheil-Small \cite{CSS} proposed the following conjecture concerning the coefficient bounds for harmonic univalent mappings.

\begin{Conj}\label{Har_Coeff}{\rm \cite[Open questions]{CSS}}
For $f=h+\overline{g} \in \mathcal{S}^0_H$ with the series representations of $h$ and $g$ as in \eqref{111}, we have
\be\label{PSCoefConj}
  |a_n| \leq \frac{(n+1)(2n+1)}{6}, ~
 |b_n| \leq \frac{(n-1)(2n-1)}{6},  ~\mbox{ and }~  \big| |a_n| - |b_n| \big| \leq n,
\ee
for all $ n \geq 2$.   The bounds are attained for the harmonic Koebe function $K(z)$, defined by
\be\label{Har_Koebe}
K(z)=\frac{z-\frac{1}{2}z^2+\frac{1}{6}z^3}{(1-z)^3}+
\overline{\left (\frac{\frac{1}{2}z^2+\frac{1}{6}z^3}{(1-z)^3}\right )}.
\ee
\end{Conj}
There are several results in support of this conjecture along with many related investigations (cf. \cite{AAP-2019,CSS,WangLZhang-01}). For example,
this conjecture has been verified for a number of subclasses of $\mathcal{S}^0_H$, namely, the classes
${\mathcal S}_H^{*0}$, ${\mathcal T}_H^0$ and ${\mathcal K}_H^0$ (hence for the class of functions convex in one direction).
 As in the case of  Bieberbach conjecture, 
Conjecture~A 
 and the recent developments on this topic
have been a driving force behind the development of univalent harmonic mappings in the plane.


\subsection{$K$-Quasiconformal harmonic mappings}
If a sense-preserving and locally univalent harmonic mapping $f=h+\overline{g}$ on $\ID$ satisfies the condition
$ |\omega_f(z)| \leq k<1$ for $\ID$,  then $f$ is called $K$-quasiregular harmonic mapping in $\ID$, where $\omega_f =g'/h'$ and
$K=\frac{1+k}{1-k}\geq 1$. We say that $f$ belongs to the class ${\mathcal S}_H(K)$ of $K$-quasiconformal harmonic mappings, if $f\in {\mathcal S}_H$ and $f$ is $K$-quasiregular harmonic mapping  in $\ID$.
 As with the standard practice, a function $f$ is called a $K$-quasiconformal harmonic mapping if it belongs to
${\mathcal S}_H(K)$ for some $K\geq 1$. Also, we define
$${\mathcal S}_H^0(K):={\mathcal S}_H(K)\cap {\mathcal S}_H^0.
$$
 For a non-trivial estimate on $|a_2|$ for $f\in {\mathcal S}_H^0(K)$, we refer to \cite[Proposition 1.6]{CP-22}. See also \cite{LP-25},
where a subclass of ${\mathcal S}_H^0(K)$ has been considered.
Recently, Wang,et al. \cite{wwrq} introduced $K$-quasiconformal harmonic Koebe function $P_k(z)$ defined by
\beq\label{eq-qhk}
P_k(z)&=&h(z)+\overline{g(z)}\\\nonumber
&=&\frac{1}{(k-1)^3}\left (\frac{(k-1)(1-3k+2kz)z}{(1-z)^2}+k(k+1)\log\Big(\frac{1-z}{1-kz}\Big)\right)\\\nonumber
&&+\frac{k}{(k-1)^3}\overline{\left(\frac{(1-k)(1+k-2z)z}{(1-z)^2}+(k+1)\log\Big(\frac{1-z}{1-kz}\Big)\right)}\\\nonumber
&=&z+\sum_{n=2}^{\infty}A(n,k)z^n+\sum_{n=2}^{\infty} B(n,k)\overline{z}^n,
\eeq
where
\beq\label{eq-mak}
A(n,k)=\frac{(1-k)n^2-2kn+k(k+1)(1+k+\cdots+k^{n-1})}{n(1-k)^2}
\eeq
and
\beq\label{eq-mbk}
B(n,k)=\frac{k(1-k)n^2-2kn+k(k+1)(1+k+\cdots+k^{n-1})}{n(1-k)^2}.
\eeq

\subsection{Image of $K$-quasiconformal harmonic Koebe function $P_k(z)$} In order to know the precise image of the unit disk under $K$-quasiconformal harmonic Koebe function $P_k(z)$, we rewrite $P_k(z)$ as
\beq \nonumber
P_k(z)&=& \frac{(1-3k+2kz)z}{(k-1)^2(1-z)^2}- \frac{k}{(k-1)^2}\overline{\left(\frac{(1+k-2z)z}{(1-z)^2}\right )} \\\nonumber
&&+\frac{k(k+1)}{(k-1)^3}\left [\log\Big(\frac{1-z}{1-kz}\Big)+ \overline{\log\Big(\frac{1-z}{1-kz}\Big)}\right]\\\nonumber
\eeq
which may be conveniently written as
$$P_k(z)=\frac{1-4k-k^2}{(1-k)^2}{\rm Re}\,\{k_e(z)\}+\frac{4k}{(1-k)^2}{\rm Re}\,\{zk_e(z)\}+2{\rm Re}\,\{l(z)\}+i{\rm Im}\,\{k_e(z)\},$$
where $k_e(z)=\frac{z}{(1-z)^2}$ is the analytic Koebe function and
$$l(z)=\frac{k(k+1)}{(k-1)^3}\log\Big(\frac{1-z}{1-kz}\Big).
$$
Note that $k_e(z)$ maps $\ID$ onto $\IC \backslash (-\infty, -1/4]$ and $k_e(e^{it}) \,(t\neq 0, 2\pi)$ belongs to $(-\infty, -1/4)$, ${\rm Im}\,\{k_e(z)\}=0$
for $|z|=1$ ($z\neq1$). So, we now try to compute $P_k(e^{it}) \,(t\neq 0, 2\pi)$. To do this, we let
$$\zeta=\zeta (z)=\frac{1+z}{1-z}=\xi+i\eta$$
which maps $\ID$ onto the right half-plane ${\rm Re}\,\{\zeta\}>0$. Calculations show that
\beqq
k_e(z)&=&\frac{\zeta ^2-1}{4}=\frac{\xi ^2-\eta ^2-1+2i\xi \eta}{4}\\
zk_e(z)&=&\frac{(\zeta -1) ^2}{4} \\
\frac{1-kz}{1-z}&=& \frac{\zeta +1}{2} -k \left (\frac{\zeta -1}{2} \right )=\frac{1+k +(1-k)\zeta}{2} \\
\left | \frac{1-kz}{1-z}\right |^2&=&  \frac{(1+k +(1-k)\xi)^2 +(1-k)^2\eta ^2}{4} \\
\eeqq
and therefore,
\beqq
P_k(z)&=&\frac{1-4k-k^2}{(1-k)^2}{\rm Re}\,\left\{\frac{\zeta^2-1}{4}\right\}+\frac{4k}{(1-k)^2}{\rm Re}\,\left\{\frac{(\zeta-1)^2}{4}\right\}\\
&&+\frac{k(k+1)}{(1-k)^3}\log\left|\frac{1+k+(1-k)\zeta}{2}\right|^2+i{\rm Im}\,\left\{\frac{\zeta^2-1}{4}\right\}\\
&=&\frac{1-4k-k^2}{(1-k)^2}\left (\frac{\xi^2-\eta^2-1}{4}\right )+\frac{4k}{(1-k)^2}\left (\frac{(\xi-1)^2-\eta^2}{4}\right )\\
&&+\frac{k(k+1)}{(1-k)^3}\log\left[\frac{\big(1+k+(1-k)\xi\big)^2}{4}+\frac{(1-k)^2\eta^2}{4}\right]+i\frac{\xi\eta}{2}.
\eeqq
If $z=e^{it}\neq1$ in $\zeta (z)=\xi+i\eta$, then $\xi=0,\, \eta\in\mathbb{R}$, and therefore, we have
\beqq
P_k(z)&=&\frac{1+k}{1-k}\left (\frac{-\eta^2}{4}\right )-\frac{1-8k-k^2}{4(1-k)^2}+\frac{k(k+1)}{(1-k)^3}\log\left[\frac{(1+k)^2+(1-k)^2\eta^2}{4}\right]\\
&\leq&\frac{k^2+8k-1}{4(1-k)^2}+\frac{2k(k+1)}{(1-k)^3}\log\frac{1+k}{2}=:M(k)\;\;(\eta=0).
\eeqq
Observe that each point $z\neq1$ on the unit circle $|z|=1$ is carried onto $(-\infty, M(k)]$.
Next if $\xi\eta=0$, i.e. $\eta=0$ (as $\xi>0$), we get
\beqq
P_k(z)&=&\frac{1-4k-k^2}{(1-k)^2}\left (\frac{\xi^2-1}{4}\right )+\frac{k(\xi-1)^2}{(1-k)^2}+\frac{2k(k+1)}{(1-k)^3}\log\left[\frac{1+k+(1-k)\xi}{2}\right]\\
&>&M(k)\;\;(\xi=0).
\eeqq
Thus under the mapping $P_k(z)$, $(-1, 1)$ is mapped monotonically onto the real interval $(M(k), +\infty)$.
Finally, each hyperbola $\xi\eta=c$, where $c\neq0$ is a real constant, is carried univalently to the set $\{w=u+i\frac{c}{2}\}$, where
\beqq
u&=&\frac{1-4k-k^2}{4(1-k)^2}\left (\xi^2-1+\frac{c^2}{\xi^2}\right )+\frac{k}{(1-k)^2}\left ((\xi-1)^2+\frac{c^2}{\xi^2}\right )\\
&&+\frac{k(k+1)}{(1-k)^3}\log\left[\frac{(1+k+(1-k)\xi)^2}{4}+\frac{(1-k)c^2}{4\xi^2}\right]
\eeqq
and $\xi>0$, which is the entire line $\{w=u+i\frac{c}{2}: -\infty<u<+\infty\}$.
This proves directly that $P_k(z)$ maps $\ID$ onto the entire plane minus the real interval $(-\infty, M(k)]$.

It is a simple exercise to see  that
$$\lim_{k\rightarrow1^-}M(k)=-\frac{1}{6}.
$$
Indeed, if we let $s=1-k$, then we have
\beqq
\lim_{k\rightarrow1^-}M(k)&=&\lim_{s\rightarrow0}\left[\Big(\frac{1}{4}-\frac{5}{2s}+\frac{2}{s^2}\Big)
+\Big(\frac{4}{s^3}-\frac{6}{s^2}+\frac{2}{s}\Big)\log\Big(1-\frac{s}{2}\Big)\right]\\
&=&\lim_{s\rightarrow0}\left[\Big(\frac{1}{4}-\frac{5}{2s}+\frac{2}{s^2}\Big)
+\Big(-\frac{2}{s^2}+\frac{5}{2s}-\frac{5}{12}-\frac{s^2}{12}+\cdots\Big)\right]\\
&=&-\frac{1}{6},
\eeqq
since $\log(1-\frac{s}{2})=-\sum_{n=1}^{\infty}\frac{s^n}{n2^n}$.

In view of the above discussion, it is natural to propose the following:

\begin{Conj}
Each function $f\in \mathcal{S}^0_H(K)$ contains the disk $\{w:\, |w|<M(k)\}$, where $M(k)=\frac{k^2+8k-1}{4(1-k)^2}+\frac{2k(k+1)}{(1-k)^3}\log\frac{1+k}{2}$.
\end{Conj}

\subsection{On the coefficient conjecture on $\mathcal{S}^0_H(K)$}
In view of the $K$-quasiconformal harmonic Koebe function $P_k(z)$, obtained by using the method of shearing by Clunie and Sheil-Small, the authors \cite{wwrq} proposed several conjectures but without much evidence in support of the validity of these conjectures. One of their conjectures reads as follows:
%
%

 \begin{Conj}\label{Har_Coeff-1}{\rm \cite{wwrq}}
For $f=h+\overline{g} \in \mathcal{S}^0_H(K)$ with the series representations of $h$ and $g$ as in \eqref{111}, we have
\be\label{PSCoefConj-1}
|a_n|\leq A(n,k),\;\;|b_n|\leq B(n,k)
~\mbox{ and }~ \big||a_n|-|b_n|\big|\leq n\;\;\mbox{for}\;\;n=2,3,\ldots,
\ee
where $A(n,k)$ and $B(n,k)$ are defined by \eqref{eq-mak} and \eqref{eq-mbk}, respectively.
Equalities occur for the function $P_k$ defined by \eqref{eq-qhk}.
\end{Conj}

Also, the authors in \cite{wwrq} remarked that it was indeed premature in view of the fact that the related conjecture for the class $\mathcal{S}^0_H$ remains an open problem. However, as with the development on the family $\mathcal{S}^0_H$, it is natural ask the following:

\bques\label{Ques1}
{\em Does Conjecture~B 
hold for every starlike function (resp. close-to-convex function) from $\mathcal{S}^0_H(K)$? How about the
validity of this conjecture for typically real $K$-quasiconformal harmonic mappings?}
\eques

It is useful to introduce the following analogous notations:
$$
\left \{\ba{rl}
{\mathcal C}_H(K)&:={\mathcal S}_H(K)\cap {\mathcal C}_H ~\mbox{ and }~{\mathcal C}_H^0(K):={\mathcal S}_H(K)\cap {\mathcal C}_H^0;\\
{\mathcal S}_H^{*}(K)&:={\mathcal S}_H(K)\cap {\mathcal S}_H^* ~\mbox{ and }~{\mathcal S}_H^{*0}(K):={\mathcal S}_H(K)\cap {\mathcal S}_H^{*0};\\
{\mathcal K}_H(K)&:={\mathcal S}_H(K)\cap {\mathcal K}_H ~\mbox{ and }~{\mathcal K}_H^{0}(K):={\mathcal S}_H(K)\cap {\mathcal K}_H^0.
\ea\right .
$$
We say that $f$ belongs to the class ${\mathcal T}_H(K)$ of typically real $K$-quasiconformal harmonic mappings, if $f\in {\mathcal T}_H$ and $f$ is a $K$-quasiconformal mapping. Moreover, we define
$${\mathcal T}_H^0(K):={\mathcal T}_H(K)\cap {\mathcal T}_H^0.$$

The main aim of this article is to give an affirmative answer to  Question \ref{Ques1}.

\begin{thm}\label{thm-main}
Let ${\mathcal F}_0$ be any one of the three classes ${\mathcal S}_H^{*0}(K)$, ${\mathcal K}_H^{0}(K)$ and ${\mathcal T}_H^0(K)$. Then Conjecture~B 
holds with  ${\mathcal F}_0$ in place of $\mathcal{S}^0_H(K)$.
\end{thm}

The constants $A(n,k)$ and $B(n,k)$ may be simplified to the following forms
$$A(n,k) =\frac{1}{n(1-k)^3}\Big(n^2+(-2n^2-2n+1)k+(n+1)^2k^2-k^{n+1}-k^{n+2}\Big)
$$
and
$$B(n,k)=\frac{k}{n(1-k)^3}\Big((n-1)^2+(-2n^2+2n+1)k+n^2k^2-k^n-k^{n+1}\Big),
$$
respectively.

In order to show that  $A(n,k)$ and $B(n,k)$ are increasing functions of $k$, we need the following lemma.

\blem\label{lem1}
For $x \in [0,1)$ and $n \geq 2$, let
$$  \phi_n(x) =\frac{M_n(x)}{(1-x)^3}
 ~\mbox{ and }~ \psi_n(x) = \frac{x}{(1-x)^3} M_{n-1}(x) ~\mbox {for all $x\in [0,1)$},
$$
where
\be\label{Mx}
 M_n(x) = n^2 + (-2n^2 - 2n + 1)x + (n+1)^2 x^2 - x^{n+1} - x^{n+2}.
\ee
Then $\phi_n$ and $\psi_n$ are strictly increasing on $[0,1)$.
\elem
\bpf First we find that
\be\label{phix}
  \phi_n'(x) = \frac{D_n(x)}{(1-x)^4} ~\mbox{for all $x\in [0,1)$},
\ee
where $D_n(x)  = (1-x)M'_n(x) + 3M_n(x)$.  In order to show that $\phi_n$ is strictly increasing on $[0,1)$, it suffices to show that $D_n(x)>0$ for all $x\in [0,1)$. To do this,
using (\ref{Mx}), we rewrite
  \begin{align*}
  	 D_n(x) & = (1-x)\big[ (-2n^2 - 2n + 1) + 2(n+1)^2x - (n+1)x^n - (n+2)x^{n+1} \big]\\
  	 		&\qquad+ 3\big[ n^2 + (-2n^2 - 2n + 1)x + (n+1)^2x^2 - x^{n+1} - x^{n+2} \big]\\
  	 		& = n^2 - 2n + 1 + x\big[ 3(-2n^2 - 2n + 1) - (-2n^2 - 2n + 1)+ 2(n+1)^2 \big] \\
  	 		& \qquad + x^2 \big[ 3(n+1)^2 - 2(n+1)^2 \big] -(1-x)x^n\big[ n+1 + (n+2)x \big] - 3x^{n+1}(1+x)\\
  	 		&
  	 		=(n-1)^2+x\big[2(-2n^2-2n+1)+2(n+1)^2\big]+x^2(n+1)^2\\
  	 		&\qquad-(n+1)x^n-4x^{n+1}+(n-1)x^{n+2}\\
  	 		&   	 		= (n-1)^2 + x(4-2n^2) + x^2(n+1)^2 - (n+1)x^n - 4x^{n+1} + (n-1)x^{n+2}.
  \end{align*}
We see that $D_n(1)= (n-1)^2>0$ for all $n\geq 2$,  and $D_n(1) = 0$. It then suffices to prove that $D_n'(x) < 0$ on $[0,1)$, so that $0=D_n(1)<D_n(x)\leq D_n(0)$ on $[0,1)$. Now, we compute
\beqq
D_n'(x) & =& 4 - 2n^2 + 2x(n+1)^2 - n(n+1)x^{n-1} - 4(n+1)x^n + (n-1)(n+2)x^{n+1},\\
D_n''(x) & =& (n+1)\big[2(n+1) - n(n-1)x^{n-2} - 4n x^{n-1} + (n-1)(n+2)x^n\big],
\eeqq
and it is easy to see that $ D_n'(0)<0$, $D_n''(0) >0$ and $ D_n'(1) =0= D_n''(1)$ for all $n\geq 2$.
Note that for $n=2$, it follows that $D_2''(x)= 12(1-x)^2$ and thus, $ D_2'''(x)= -24(1-x) < 0$ for all $x\in [0,1)$. Moreover, for $n \geq 3$, we see that
\begin{align*}
	D_n'''(x) & = (n+1)\big[ -n(n-1)(n-2)x^{n-3} - 4n(n-1)x^{n-2} + (n-1)n(n+2)x^{n-1} \big]\\
	& = (n-1)n(n+1)x^{n-3}\big[-(n-2) - 4x + (n+2)x^2 \big]\\
	&= -(n-1)n(n+1)x^{n-3}(1-x)\big[n-2 + (n+2)x\big]\\
	&\leq 0\;\;\mbox{for all $x\in [0,1)$}.
\end{align*}
Thus, we obtain that for all $n \geq 2$, we have
$$D_n'''(x) \leq 0 ~\mbox{ for all $x\in [0,1)$}
$$
which gives that
$$ D_n''(x) > D_n''(1) = 0 ~\mbox {for all $x\in [0,1)$}
$$
Therefore, $D_n'(x) < D_n'(1) = 0$ for all $x\in [0,1)$ and this in turn shows that
$D_n(x) > D_n(1) = 0$ for all $x\in [0,1)$. Hence, by (\ref{phix}), the function $\phi_n(x)$ is increasing (strictly) increasing on $[0,1)$.

Next, we consider the second part of the lemma. We now compute
$$
\psi_n'(x) = \frac{E_n(x)}{(1-x)^4} ~\mbox {for all $x\in [0,1)$},
$$
where
$$
M_{n-1}(x) = (n-1)^2 + (-2n^2+2n+1)x + n^2x^2 - x^n - x^{n+1},
$$
and
$$
E_n(x) = (1-x)\big[x M_{n-1}'(x) + M_{n-1}(x)\big] + 3x M_{n-1}(x).
$$
To show that $\psi_n$ is strictly increasing on $[0,1)$, it suffices to show that $E_n(x)> 0$ on $[0,1]$. Now, we may rewrite $E_n(x)$ as
$$
E_n(x) = x\big[(1-x) M_{n-1}'(x) + 3M_{n-1}(x)\big] + (1-x)M_{n-1}(x)
$$
and thus,
\begin{equation*}
 	E_n(x)= xD_{n-1}(x) + (1-x)M_{n-1}(x).
\end{equation*}
First, we note that for $n=2$,
\begin{align*}
E_2(x)& = xD_1(x) + (1-x)M_1(x)\\
& = x\big[(1-x)M_1'(x) + 3M_1(x)\big] + (1-x)M_1(x)\\
& = x\big[(1-x)(-3+6x-3x^2) + 3(1-3x+3x^2-x^3)\big] \\
&\qquad + (1-x)(1-3x+3x^2+x^3)\\
 		&= x\big[-3(1-x)(1-x)^2 + 3(1-x)(1-2x+x^2)\big] + (1-x)(1-x)^3\\
 		&= (1-x)^4 > 0 ~\mbox{ for all $x\in [0,1)$.}
\end{align*}
Next, we consider the case $n \geq 3$ and find that
\begin{align*}
M_{n-1}'(x) &=(-2n^2+2n+1) + 2n^2x - n x^{n-1} - (n+1)x^n,\\
M_{n-1}''(x) &= 2n^2 - n(n-1)x^{n-2} - n(n+1)x^{n-1},
\end{align*}
and thus, it is clear that $M_{n-1}'''(x) < 0$ for all $x \in [0,1)$ and $n\geq 3$, showing that
$$M_{n-1}''(x) > M_{n-1}''(1)=0 ~\mbox{ for all $x\in [0,1)$},
$$
which implies that
$$M_{n-1}'(x) < M_{n-1}'(1) = 0 ~\mbox{ for all $x\in [0,1)$.}
$$
Thus,
$$M_{n-1}(x) > M_{n-1}(1) = 0 ~\mbox{ for all $x\in [0,1)$.}
$$
The last observation together with the fact that $D_n(x)> 0$ on $[0,1)$, shows that $\psi_n'(x)> 0$ on $[0,1)$. This completes the proof of the lemma.
\epf

\begin{rem}\label{rem-3}
By Lemma \ref{lem1},  $\phi_n(k)=nA(n,k)$ and $\psi_n(k)=nB(n,k)$ and thus, $A(n,k)$ and $B(n,k)$ are increasing functions of $k\in [0,1)$. If $k\rightarrow1^-$, then Theorem \ref{thm-main} shows that
$$\lim_{k\rightarrow1^-}A(n,k)=\frac{(2n+1)(n+1)}{6}
~\mbox{ and }~
\lim_{k\rightarrow1^-}B(n,k)=\frac{(2n-1)(n-1)}{6}.$$
In this limiting case, Theorem \ref{thm-main} clearly coincides with Conjecture~A 
and the results for three geometric subfamilies ${\mathcal S}_H^{*0}$,  ${\mathcal K}_H^0$ and ${\mathcal T}_H^0$ obtained in \cite{DuR}.

The case $n=2$ of Theorem \ref{thm-main} gives the sharp inequalities
$$|a_2|\leq \frac{5K+3}{2K+2}\;\;\mbox{and}\;\;|b_2|\leq \frac{K-1}{2K+2}.
$$
Hence, Theorem \ref{thm-main} also solves Conjecture~B 
of \cite{wwrq} for the subfamilies ${\mathcal S}_H^{*0}$,  ${\mathcal K}_H^0$ and ${\mathcal T}_H^0$, which include  class of
typically real $K$-quasiconformal harmonic mappings $f$ which do not require the univalency assumption on $f$.
\end{rem}

\bcor \label{cor-main}
Let ${\mathcal F}$ be any one of the  three classes ${\mathcal S}_H^{*}(K)$, ${\mathcal K}_H(K)$ and ${\mathcal T}_H(K)$.
Then the coefficients of each function $f\in {\mathcal F}$ $(K>1)$ satisfy the sharp inequalities
$$|a_n|\leq E_n(k):= A(n, k_0) + \frac{1-\sqrt{1-k_0^2}}{k_0}  B(n, k_0),$$
$$|b_n|\leq F_n(k):=B(n, k_0) + \frac{1-\sqrt{1-k_0^2}}{k_0} A(n, k_0),
$$
for $n=2,3,\ldots,$ where $k_0=\frac{2k}{1+k^2}$, $A(n,k)$ and $B(n,k)$ are defined by \eqref{eq-mak} and \eqref{eq-mbk}, respectively.
Equalities occur for the function $Q_k$ defined by
\beq\label{eq-qhkq}
Q_k(z)=P_k(z)+\frac{\sqrt{K_0}-1}{\sqrt{K_0}+1}\overline{P_k(z)},\eeq
where $P_k(z)\in {\mathcal F}_0^0$ is defined by \eqref{eq-qhk}, ${\mathcal F}_0^0$ is any one of the  three classes ${\mathcal S}_H^{*0}(K_0)$, ${\mathcal K}_H^0(K_0)$ and ${\mathcal T}_H^0(K_0)$ and $K_0=K^2$.
\ecor

The limiting case of Corollary \ref{cor-main} as $k\rightarrow1^-$ gives the following well-known result.

\bcor \label{cor-main2}
Let ${\mathcal F}$ be any one of the three classes ${\mathcal S}_H^{*}$, ${\mathcal K}_H$ and ${\mathcal T}_H$.
Then the coefficients of each function $f=h+\overline{g}\in {\mathcal F}$   satisfy the sharp inequalities
$$|a_n|< \frac{2n^2+1}{3}  ~\mbox{ and }~ |b_n|< \frac{2n^2+1}{3} ~\mbox{ for $n=2,3,\ldots.$ }~
$$
\ecor
\bpf
As remarked earlier, Lemma \ref{lem1} implies that  $A(n,k_0)$ and $B(n,k_0)$ are increasing functions of $k_0\in [0,1)$.
Thus, both
$$A(n, k_0) + \frac{1-\sqrt{1-k_0^2}}{k_0}  B(n, k_0),
~\mbox{ and }~
B(n, k_0) + \frac{1-\sqrt{1-k_0^2}}{k_0} A(n, k_0),
$$
are increasing functions of $k_0\in (0,1)$. Also, we observe that
$$ k_0\mapsto \frac{1-\sqrt{1-k_0^2}}{k_0}
$$
is increasing on  $k_0\in (0,1)$, $k_0\ra 0$ if and only if $k\ra 0$, and $k\ra 1$ if and only if $k_0\ra 1$.
Hence,
$$\lim_{k\rightarrow1^-}  E_n(k) = \frac{2n^2+1}{3}=\lim_{k\rightarrow1^-}  F_n(k)
$$
and the proof follows.
\epf

In order to prove our main result, we first need to consider the natural analogue of Conjecture~B 
for $f  \in \mathcal{S}^0_H(K)$ that
are convex in $\ID$. Secondly we  separate the proof of the typically real case in the form of a separate theorem,
as functions in ${\mathcal T}_H(K)$ and ${\mathcal T}_H^0(K)$ are not necessarily univalent in $\ID$.

\br
Area minimizing property of ${\mathcal S}_H^0(K)$ follows easily as with the case of ${\mathcal S}_H^0$. Indeed for each $f\in {\mathcal S}_H^0(K)$ area of $f(\ID)\geq \pi \Big(1-\frac{k^2}{2}\Big)$ and
the minimum is attained only for $f_0(z)=z+\frac{k}{2}{\overline z}^2$ and its rotations, where $k=(K-1)/(K+1)$. For the sake of completeness and clarity, we
simply brief the proof here (cf. \cite{Bo}). It suffices to recall that
the area $A$ of $f(\ID)$  is given by
$$ A =\iint_{\ID} J_f(z)\,dx\,dy  \geq  \iint_{\ID}\big(|h'(z)|^2-k^2|zh'(z)|^2\big) \,dx\,dy=:B,
$$
%
where
$$B=\pi\sum_{n=1}^\infty n\Big(1-\frac{k^2n}{n+1}\Big)|a_n|^2\geq   \pi \Big(1-\frac{k^2}{2}\Big).
$$
 The minimum is attained if we choose $a_n=0$ for all $n\geq 2$.
This happens when $\omega(z)=e^{i\alpha}kz$ for some $\alpha$, which gives
$g'(z)=ke^{i\alpha}zh'(z)=ke^{i\alpha}z$ so that
$$f(z)=z+\frac{k}{2}e^{-i\alpha}\overline{z}^2.
$$
\er

\subsection{Convex $K$-quasiconformal harmonic mappings}
Throughout the discussion, $a_n:=a_n(h)$ and $b_n:=b_n(g)$ denote respectively the Taylor coefficients of $h$ and $g$ in the representation of $f=h+\overline{g}\in \mathcal{H}$.
%
We begin the discussion by recalling the following results due to Clunie and Sheil-Small \cite{CSS}.

\begin{Thm}\label{thmA}
The coefficients of each function $f\in {\mathcal C}_H^0$ satisfy the sharp inequalities
$$|a_n|\leq \frac{n+1}{2},\;\;|b_n|\leq \frac{n-1}{2}\;\;\mbox{and}\;\;\big ||a_n|-|b_n|\big|\leq 1$$
for $n=2,3,\ldots.$ Equality occurs for the function $L$, where
\beqq
L(z)&=&\frac{1}{2}\Big(\frac{z}{1-z}+\frac{z}{(1-z)^2}\Big)+\overline{\frac{1}{2}\Big(\frac{z}{1-z}-\frac{z}{(1-z)^2}\Big)}\\
&=&z+\sum_{n=2}^{\infty}\frac{n+1}{2}z^n-\sum_{n=2}^{\infty}\frac{n-1}{2}\overline{z}^n.
\eeqq
Also, the coefficients of each function $f\in {\mathcal C}_H$ satisfy the sharp inequalities
$$|a_n|<n\;\;\mbox{and}\;\;  |b_n|<n ~\mbox{ for $n=2,3,\ldots.$}
$$
\end{Thm}

%


Now, we consider the coefficients estimate of $f\in {\mathcal C}_H(K)$ and also for $f\in {\mathcal C}_H^0(K)$.

\begin{thm}\label{thm1}
The coefficients of each function $f\in {\mathcal C}_H^0(K)$ satisfy the sharp inequalities
$$|a_n|\leq a(n,k),
\;\; |b_n|\leq b(n,k) ~\mbox{ and }~ \big||a_n|-|b_n|\big|\leq 1\;\;\mbox{for}\;\;n=2,3,\ldots,
$$
where
$$a(n,k)=\frac{n-k(n+1)+k^{n+1}}{n(1-k)^2}~\mbox{ and }~b(n,k)=\frac{k}{n(1-k)^2}(n-1-nk+k^n).
$$
Specially, if $n=2$, we have the sharp inequalities
$$|a_2|\leq  \frac{3K+1}{2K+2}\;\;\mbox{and}\;\;|b_2|\leq \frac{K-1}{2K+2}.$$
Equalities occur for the function $P\in {\mathcal C}_H^0(K)$ defined by
\beq\label{eq2}
P(z)=h_k(z)+\overline{g_k(z)}=z+\sum_{n=2}^{\infty}a(n,k)z^n-\sum_{n=2}^{\infty}b(n,k)\overline{z}^n,
\eeq
where
$$h_k(z)=\frac{z}{(1-k)(1-z)}-\frac{k}{(1-k)^2}\log\left(\frac{1-kz}{1-z}\right),
$$
and
$$g_k(z)=
\frac{-kz}{(1-k)(1-z)}-\frac{k}{(1-k)^2}\log\left(\frac{1-z}{1-kz}\right).
$$

\end{thm}

Theorem \ref{thm1} leads also to sharp coefficient bounds for mappings of the class ${\mathcal C}_H(K)$. In order to state this result, we need the following lemma.

\blem\label{lem2}
For $x \in [0,1)$ and $n \geq 2$, let
$$  G_n(x) =\frac{L_n(x)}{(1-x)^2}
 ~\mbox{ and }~ H_n(x) = \frac{x}{(1-x)^2} L_{n-1}(x) ~\mbox {for all $x\in [0,1)$},
$$
where
\be\label{Lx}
 L_n(x) = n - (n + 1)x + x^{n+1}
\ee
Then $G_n$ and $H_n$ are strictly increasing on $[0,1)$.
\elem
\bpf
Clearly, $L_n(x) =n(1-x)-x(1-x^n)=(1-x)[n-x(1+x+\cdots x^{n-1})]>0$ on $[0,1)$. Moreover,
$$L_n'(x) =-(n+1)(1-x^n), ~L_n''(x) =n(n+1)x^{n-1}  ~\mbox{ and }~L_n'''(x) =(n-1)n(n+1)x^{n-2}
$$
so that $L_n'(x)<0$ on $[0,1)$, $L_n''(x)>0$ on $(0,1)$, and $L_n'''(x)>0$ on $(0,1)$. Moreover, $L_n'(1)=0$ and $L_n''(1)>0$. Now, a computation gives
$$
G_n'(x) = \frac{F_n(x)}{(1-x)^3} ~\mbox{for all $x\in [0,1)$},
$$
where $F_n(x)  = (1-x)L'_n(x) + 2L_n(x)$. Again,
$$F_n'(x)  = (1-x)L''_n(x) + L'_n(x) ~\mbox{ and }~ F_n''(x) = (1-x)L'''_n(x),
$$
showing that $F_n''(x)\geq 0$ for all $x\in [0,1)$ and $n\geq 2$. Therefore, $F_n'(x)\leq F_n'(1)=0$  for all $x\in [0,1)$, which in turn implies that
$F_n(x)\geq F_n(1)=0$  for all $x\in [0,1)$, Thus, $G_n(x)$ is strictly increasing on $(0,1)$.

For the proof the second part, we find that $H_2(x)=x$ and so, there is nothing to prove. We may thus, assume that $n\geq 3$ and we have
$$
H_n'(x) = \frac{J_n(x)}{(1-x)^3} ~\mbox{for all $x\in [0,1)$},
$$
where $J_n(x)  = (1-x)[xL'_{n-1}(x) + L_{n-1}(x)] +2xL_{n-1}(x)$. We need to show that  $J_n(x)\geq 0$ on $[0,1)$. Indeed,  $J_n(0)=n-1>0$ and we may rewrite  $J_n(x)$ as
$$J_n(x)  = x[(1-x)L'_{n-1}(x) + 2L_{n-1}(x)] +(1-x)L_{n-1}(x).
$$
By the previous argument, it is clear that $J_n(x)>0$ on $(0,1)$ and for $n\geq 3$. Thus, $H_n(x)$ is strictly increasing on $(0,1)$.
\epf

\begin{rem}\label{rem-1}
First we observe  that $nG_n(k)= a(n,k)$ and $nH_n(k)= b(n,k)$. By Lemma \ref{lem2}, $a(n,k)$ and $b(n,k)$  are increasing functions of $k\in [0,1)$.
If $k\rightarrow1^-$, then
$$\lim_{k\rightarrow1^-} a(n,k)=\frac{n+1}{2} ~\mbox{ and  }~ \lim_{k\rightarrow1^-} b(n,k)=\frac{n-1}{2}.
$$
Hence, in this case, Theorem \ref{thm1} coincides with the first part of Theorem~C. 
\end{rem}

\begin{cor}\label{cor1}
The coefficients of each function $f\in {\mathcal C}_H(K)$ $(K>1)$ satisfy the sharp inequalities
$$|a_n|\leq C_n(k) = a(n,k_0) + \left (\frac{1-\sqrt{1-k_0^2}}{k_0} \right )b(n,k_0)
$$
and
$$|b_n|\leq D_n(k) = b(n,k_0) + \left ( \frac{1-\sqrt{1-k_0^2}}{k_0} \right ) a(n,k_0)
$$
for $n=2,3,\ldots$, where $k_0=\frac{2k}{1+k^2}$.  Specially, if $n=2$, we have the sharp inequalities
$$|a_2|\leq C_2(k)= \frac{1+k+2k^2}{1+k^2}=\frac{2K^2-K+1}{K^2+1}
$$
and
$$ |b_2|\leq D_2(k)=\frac{k(2+k+k^2)}{1+k^2}= \frac{(K-1)(2K^2+K+1)}{(K+1)(K^2+1)}.
$$
Equalities occur for the function $Q$ defined by
\be\label{eq3}
Q(z)=P(z)+\frac{\sqrt{K_0}-1}{\sqrt{K_0}+1}\overline{P(z)},
\ee
where $P(z)\in {\mathcal C}_H^0(K_0)$ is defined by \eqref{eq2} and $K_0=K^2$.
\end{cor}

\begin{rem}\label{rem-2}
As before, we remark that $C_n(k)$ and $D_n(k)$ are strictly increasing functions of $k\in [0,1)$.
If $k\rightarrow1^-$, then, as in the case of the proof of Corollary \ref{cor-main2}, we have
$$\lim_{k \rightarrow1^-} C_n(k)=n ~\mbox{ and } ~ \lim_{k\rightarrow1^-} D_n(k)=n.
$$
%
Thus,  in the limiting case, Corollary \ref{cor1} coincides with the second part of Theorem~C. 
\end{rem}

\subsection{Typically real $K$-quasiconformal harmonic mappings}

By \eqref{eq-qhk}, since $P_k(0)=0$, $h'(0)+\overline{g'(0)}=1>0$ and $\frac{g'(z)}{h'(z)}=kz$, we know that $P_k\in {\mathcal T}_H^0(K)$ (cf. \cite{BHH} or \cite[Theorem A]{WLZ}). We consider the coefficients estimate of $f\in {\mathcal T}_H(K)$ and $f\in {\mathcal T}_H^0(K)$ and show that Conjecture~B 
 holds for ${\mathcal T}_H^0(K)$. More precisely, we have

\begin{thm}\label{thm2}
If $f=h+\overline{g}\in {\mathcal T}_H^0(K)$, then $a_1=1$, and we have the sharp inequalities stated in \eqref{PSCoefConj-1} hold.
\end{thm}

%



We  present the proofs of  Theorems \ref{thm-main}, \ref{thm1}, and \ref{thm2}
in Section \ref{sec-2}.

\section{Proofs of the main results}\label{sec-2}

First, we construct the following $K$-quasiconformal sense-preserving harmonic function
that is convex in the vertical direction.
This construction is based on  the method of shearing introduced by Clunie and Sheil-Small \cite{CSS}. See also \cite{DuR,PonRasi2013}. Let
$$h(z)+g(z)=\frac{z}{1-z}\;\;\mbox{and}\;\;\omega(z)=\frac{g'(z)}{h'(z)}=ke^{i\alpha}z,\;\;\alpha\in[0,2\pi), k\in[0, 1).$$
Differentiation of the first equation together with the second relation gives the pair of linear differential equations
$$g'(z)+h'(z)=\frac{1}{(1-z)^2} ~\mbox{ and }~ke^{i\alpha}zh'(z)-g'(z)=0.
$$
After solving these two equations we get the unique solution $P_{k}^{\alpha}(z)=h(z)+\overline{g(z)}$, where
$$h'(z)=\frac{1}{(1-z)^2(1+ke^{i\alpha}z)}~\mbox{ and }~ g'(z)=\frac{ke^{i\alpha}z}{(1-z)^2(1+ke^{i\alpha}z)}.
$$
Integration now produces the expressions for $h$ and $g$:

$$h(z)=l_1(z)+l_2(z)
~\mbox{ and }~ g(z)=ke^{i\alpha}l_1(z)-l_2(z)
$$
where
$$l_1(z)=\frac{z}{(1+ke^{i\alpha})(1-z)}~\mbox{ and }~l_2(z)=\frac{ke^{i\alpha}}{(1+ke^{i\alpha})^2}\log\left(\frac{1+ke^{i\alpha}z}{1-z}\right).
$$
Then
\beqq
P_{k}^{\alpha}(z)&=& {\rm Re}\, (h(z)+g(z)) +i {\rm Im}\, (h(z)-g(z))\\
&=& {\rm Re}\, \{(1+ke^{i\alpha})l_1(z)\}+i{\rm Im}\, \{(1-ke^{i\alpha})l_1(z)+2l_2(z)\}.
\eeqq
Now, if we let $e^{i\alpha}=-1$, then a simple computation shows that
\beqq
\{P_{k}^{\pi}(z) &=& {\rm Re}\left ( \frac{z}{1-z} \right )+i \left [  \frac{1+k}{1-k} {\rm Im}\left ( \frac{z}{1-z} \right )-  \frac{2k}{(1-k)^2}
  {\rm Im}\left (  \log\left(\frac{1-kz}{1-z}\right)\right ) \right ]
\eeqq
from which it is clear that $P_{k}^{\pi}(z)$ maps $\ID$ onto the right half-plane ${\rm Re}\,w>-\frac{1}{2}$, since $\frac{z}{1-z}$ maps $\ID$ onto the right half-plane ${\rm Re}\,w>-\frac{1}{2}$. Hence, $P_{k}^{\pi}(z)$ is a univalent right
half-plane mapping for each $k\in [0,1)$ and thus, the function $P_{k}^{\pi}$ belongs to ${\mathcal C}_H^0(K)$.

Finally, a computation gives
\beq\label{eq1}
P_{k}^{\alpha}(z)&=& h(z)+\overline{g(z)}
=z+\sum_{n=2}^{\infty}a(n,k,\alpha)z^n-\sum_{n=2}^{\infty}\overline{b(n,k,\alpha)}\overline{z}^n,
\eeq
where
$$
\left \{ \begin{array}{l}
\ds a(n,k,\alpha)=\frac{n+(n+1)ke^{i\alpha}-(-1)^nk^{n+1}e^{i(n+1)\alpha}}{n(1+ke^{i\alpha})^2}, ~\mbox{ and }\\
\ds b(n,k,\alpha)=-\frac{(n-1)ke^{i\alpha}+nk^2e^{2i\alpha}+(-1)^nk^{n+1}e^{i(n+1)\alpha}}{n(1+ke^{i\alpha})^2}.
\end{array} \right.
$$

In order to present the proof of Theorem \ref{thm1}, let us recall some useful results.

\begin{Lem}$($\cite[p.~50]{DuR}$)$\label{lemA}
If $f=h+\overline{g}\in {\mathcal C}_H$, then there exist angles $\alpha$ and $\beta$ such that
$${\rm Re}\, \{(e^{i\alpha}h'(z)+e^{-i\alpha}g'(z))(e^{i\beta}-e^{-i\beta}z^2)\}>0 ~\mbox{ for all $z\in\ID$.}
$$
\end{Lem}

\begin{Lem}$($\cite[p.~51]{DuR}$)$\label{lemB}
If $\varphi(z)=c_0+c_1z+\cdots$ is analytic with ${\rm Re}\,\{\varphi(z)\}>0$ in $\ID$, then $|c_n|\leq2{\rm Re}\,\{c_0\}$, $n=1,2,\ldots$.
\end{Lem}

%



\subsection{Proof of  Theorem \ref{thm1}}
For $z\in\ID$, we see that
$$\frac{1+z}{1-z}=1+\sum_{n =1}^{\infty} 2 z^n.
$$
Then, it follows from Lemmas~D and E 
 that the Taylor coefficients of the function
$$F(z)=(e^{i\alpha}h'(z)+e^{-i\alpha}g'(z))(e^{i\beta}-e^{-i\beta}z^2)$$
are dominated in modulus by the corresponding coefficients of the function $\frac{1+z}{1-z}$.
Hence, the coefficients of $e^{i\alpha}h'(z)+e^{-i\alpha}g'(z)$ are dominated in modulus by those of
$$\frac{1+z}{1-z}\frac{1}{1-z^2}=\frac{1}{(1-z)^2}.$$
Integration shows that the coefficients of $e^{i\alpha}h(z)+e^{-i\alpha}g(z)$ are dominated in modulus by those of
$z/(1-z)=\sum_{n =1}^{\infty} z^n.$ Thus,
$$\big||a_n|-|b_n|\big|\leq |e^{i\alpha}a_n+e^{-i\alpha}b_n|\leq 1,$$
as the theorem asserts.

To obtain the other estimates, we note that the dilatation $\omega(z)=g'(z)/h'(z)$ satisfies $\omega(0)=0$ and $|\omega(z)|\leq k$ and
thus, by the classical Schwarz Lemma,  we have
$|\omega(z)|\leq k|z|$ for $z\in \ID$ and $|\omega'(0)|\leq k$. It follows that
$$\frac{\omega(z)}{e^{i\alpha}+e^{-i\alpha}\omega(z)}\prec \frac{kz}{e^{i\alpha} +ke^{-i\alpha}z},
$$
where $\prec $ denotes the usual subordination  (cf. \cite{CSS,Du,DuR}). Writing
$$g'(z)=\frac{\omega(z)}{e^{i\alpha}+e^{-i\alpha}\omega(z)}\frac{1}{e^{i\beta}-e^{-i\beta}z^2}F(z)$$
and appealing to Schwarz Lemma (applied to $\omega (z)/k$), 
we get that the coefficients of $g'(z)$ are dominated in modulus by those of
$$\frac{kz}{1-kz}\frac{1}{1-z^2}\frac{1+z}{1-z}=\frac{kz}{1-kz}\frac{1}{(1-z)^2}
=\sum_{n =0}^{\infty}\left (\sum_{m =0}^{n}k^{m+1}(n-m+1)\right ) z^{n+1}.
$$
Thus, as $b_1=0$ and $g'(z)=\sum_{n =0}^{\infty}(n+2)b_{n+2}z^{n+1}$, the above observation yields that
$$(n+2)|b_{n+2}|\leq\sum_{m =0}^{n}k^{m+1}(n-m+1) ~\mbox{ for $n\geq 0$,}
$$
that is,
$$|b_n|\leq \frac{1}{n}\sum_{m =1}^{n-1}k^{m}(n-m)=\sum_{m =1}^{n-1}k^{m}-\frac{1}{n}\sum_{m =1}^{n-1}mk^{m} ~\mbox{ for $n\geq2$.}
$$
It follows from
\beq\label{eq-k}
\sum_{m =1}^{n-1}k^{m}=\frac{k(1-k^{n-1})}{1-k}
\eeq
and
\beq\label{eq-mk}
\sum_{m =1}^{n-1}mk^{m}=\frac{k}{(1-k)^2}\Big(1-nk^{n-1}+(n-1)k^n\Big)
\eeq
that
$$|b_n|\leq \frac{k}{n(1-k)^2}(n-1-nk+k^n).$$
Then
$$
|a_n|\leq\big||a_n|-|b_n|\big|+|b_n|\leq \frac{n-(n+1)k+k^{n+1}}{n(1-k)^2}.
$$

Hence, if $n=2$, we have
$$|a_2|\leq \frac{k+2}{2}=\frac{3K+1}{2K+2}\;\;\mbox{and}\;\;|b_2|\leq \frac{k}{2}=\frac{K-1}{2K+2}.$$
To see that the bounds are sharp, we let $e^{i\alpha}=-1$ in \eqref{eq1}. Then, we obtain the function $P(z)\in {\mathcal C}_H^0(K)$
defined by \eqref{eq2}, namely,
\be\label{eq3-P}
P(z)=P_{k}^{\pi}(z)=h_k(z)+\overline{g_k(z)}=z+\sum_{n=2}^{\infty}a(n,k)z^n-\sum_{n=2}^{\infty}b(n,k)\overline{z}^n,
\ee
where $a(n,k)$ and $b(n,k)$ are given as in the statement. This completes the proof.
\qed
\medskip

\subsection{Proof of Corollary \ref{cor1}}
Let $f=h+\overline{g} \in {\mathcal C}_H(K)$,  where $|\omega _{f}(z)|\leq k$ and $k=\frac{K-1}{K+1}$. In particular, $|\omega _{f}(0)|=|b_1|\leq k$. Clearly, each $f=h+\overline{g} \in {\mathcal C}_H(K)$ has the form
$$f=f_0+\overline{b_1f_0}=h_0+ \overline{b_1}g_0+\overline{g_0 + b_1h_0}, ~ \mbox{i.e., }~ f_0=\frac{1}{1-|b_1|^2}\left [ h-\overline{b_1} g +\overline{g-b_1h}\right ],
$$
where $f_0\in {\mathcal C}_H^0(K_0)$, $|\omega _{f_0}(z)|\leq k_0$ for all $z\in \ID$ and for some $k_0<1$, and $k_0=\frac{K_0-1}{K_0+1}$. To obtain $k_0$, we use above relation and obtain that
$$\omega _f(z)=\frac{g'(z)}{h'(z)}  =
\frac{\frac{g_0'(z)}{h_0'(z)}+b_1}{ 1+\overline{b_1}\frac{g_0'(z)}{h_0'(z)}} =\frac{\omega _{f_0}(z) +b_1}{ 1+\overline{b_1}\omega _{f_0}(z)} , ~ \mbox{i.e., }~
\omega _{f_0}(z)=\frac{\omega _{f}(z) -b_1}{ 1-\overline{b_1}\omega _{f}(z)},
$$
and therefore,
$$|\omega _{f_0}(z)|\leq \frac{|\omega _{f}(z)|+ |b_1|}{ 1+|b_1|\,|\omega _{f}(z)|}\leq \frac{2k}{1+k^2}=:k_0.
$$
Here we used the fact that $x\mapsto \frac{k+x}{1+kx}$ is an increasing function of $x\in [0,k]$. Note that $k\ra 1$ if and only if $k_0\ra 1$. Moreover, solving the equation
$$\frac{2k}{1+k^2}=k_0$$
for $k$ gives
$$ k=\frac{1-\sqrt{1-k_0^2}}{k_0},
$$
since the other root is bigger than $1$. Furthermore,  $f=f_0+\overline{b_1f_0}$ gives
$$ a_n(h)=a_n(h_0) +\overline{b_1}b_n(g_0) ~\mbox{ and }~b_n(g)=b_n(g_0) +b_1a_n(h_0).
$$
Finally, Theorem \ref{thm1} and the fact that $|b_1|\leq k=\frac{1-\sqrt{1-k_0^2}}{k_0}$ give

\beqq
|a_n(h)| &\leq&\frac{n-(n+1)k_0+k_0^{n+1}}{n(1-k_0)^2}+|b_1|\frac{k_0}{n(1-k_0)^2}(n-1-nk_0+k_0^n)\\
&= & a(n,k_0) + |b_1| b(n,k_0)\\
&\leq & a(n,k_0) + \left (\frac{1-\sqrt{1-k_0^2}}{k_0}\right ) b(n,k_0) =:C_n(k)\\
\eeqq
and
\beqq
|b_n(g)|&\leq&\frac{k_0}{n(1-k_0)^2}(n-1-nk_0+k_0^n)+|b_1|\frac{n-(n+1)k_0+k_0^{n+1}}{n(1-k_0)^2}\\
&= & b(n,k_0) + |b_1| a(n,k_0)\\
&\leq & b(n,k_0) + \left ( \frac{1-\sqrt{1-k_0^2}}{k_0}\right) a(n,k_0)=:D_n(k),
\eeqq
for $n=2, 3, \ldots$, where  $a(n,k_0)$ and $b(n,k_0)$ are as in Theorem \ref{thm1}, and $ k_0=\frac{2k}{1+k^2}$.

By Lemma \ref{lem2}, $a(n,k_0)$ and $b(n,k_0)$ are increasing functions of $k_0\in [0,1)$.  Also, it is easy to observe that $k_0\mapsto \frac{1-\sqrt{1-k_0^2}}{k_0}$ is increasing for $k_0\in [0,1)$,
and $k\mapsto \frac{2k}{1+k^2}$ is also increasing for $k\in [0,1)$. Consequently, both $C_n(k)$ and $D_n(k)$ are increasing functions of $k\in [0,1)$.

 In particular, by a simple computation, the case $n=2$ of the first inequality above gives
$$|a_2(h)|\leq \frac{2+k_0}{2} +\frac{1-\sqrt{1-k_0^2}}{2}=\frac{1+k+2k^2}{1+k^2}=\frac{2K^2-K+1}{K^2+1}, 
$$
since
$$1+k+2k^2=\frac{2(2K^2-K+1)}{(K+1)^2}~\mbox{ and }~1+k^2=\frac{2(K^2+1)}{(K+1)^2}.
$$
Similarly, it is easy to see that
$$
|b_2(g)|\leq \frac{k_0}{2} +\frac{1-\sqrt{1-k_0^2}}{k_0}\left ( \frac{2+k_0}{2}\right )= \frac{k(2+k+k^2)}{1+k^2}= \frac{(K-1)(2K^2+K+1)}{(K+1)(K^2+1)},
$$
since
$$ k(2+k+k^2)=\frac{2(K-1)(2K^2+K+1)}{(K+1)^3}.
$$

To see that the bounds are sharp, we only need consider functions $Q(z)$ given by \eqref{eq3}, since
$$\frac{2k}{1+k^2}=k_0~\Leftrightarrow~K^2=K_0,$$
and the last two way implications is easy to verify.
%
%
\qed

\subsection{Proof of  Theorem \ref{thm2}}
Since $b_1=0$, the condition $f(r)>0$ for $0<r<1$ shows that $a_1>0$ and so, $a_1=1$. Consequently, $\varphi=h-g\in {\mathcal T}$ and so (see \cite[Corollary, p.~58]{Du})
\be\label{eq-extra1}
\big||a_n|-|b_n|\big|\leq |a_n-b_n|\leq n.
\ee
Again, we observe that the dilatation $\omega=g'/h'$ satisfies $|\omega(z)|\leq k|z|$ for $z\in\ID$. Now, since
$$g'(z)=\omega(z)h'(z)=\omega(z)\big(\varphi'(z)+g'(z)\big),
$$
we have
$$g'(z)=\frac{\omega(z)}{1-\omega(z)}\varphi'(z).
$$
As in the previous argument, it follows that  the coefficients of $g'(z)$ are dominated in modulus by those of the function
$$\frac{kz}{1-kz}\Big(\frac{z}{(1-z)^2}\Big)'=\sum_{n =0}^{\infty}\left (\sum_{m =0}^{n}k^{m+1}(n-m+1)^2\right ) z^{n+1}
$$
from which we easily obtain that
$$|(n+2)b_{n+2}|\leq\sum_{m =0}^{n}k^{m+1}(n-m+1)^2=\sum_{m =1}^{n+1}k^{m}(n-m+2)^2 ~\mbox{ for $n\geq 0$},
$$
or equivalently,
$$|b_{n}|\leq\frac{1}{n}\sum_{m =1}^{n-1}k^{m}(n-m)^2 ~\mbox{ for $n\geq 2$}.
$$
By simple computations, it follows from \eqref{eq-k}, \eqref{eq-mk} and
$$\sum_{m =1}^{n-1}m^2k^{m}=\frac{k}{(1-k)^3}\Big(1+k-n^2k^{n-1}+(2n^2-2n-1)k^n-(n-1)^2k^{n+1}\Big)
$$
that
$$|b_{n}|\leq\frac{k}{n(1-k)^3}\Big((n-1)^2+(-2n^2+2n+1)k+n^2k^2-k^n-k^{n+1}\Big)=B(n,k).
$$
The final estimate is easy to obtain by combining the last coefficients inequality on $|b_{n}|$ and \eqref{eq-extra1}:
\beqq
|a_{n}|&\leq& \big||a_n|-|b_n|\big|+|b_{n}|\\
&\leq& n+ B(n,k)\\
&=& \frac{1}{n(1-k)^3}\Big(n^2+(-2n^2-2n+1)k+(n+1)^2k^2-k^{n+1}-k^{n+2}\Big)\\
&=&A(n,k).
\eeqq
This completes the proof.
\qed
\medskip

\subsection{Proof of Corollary \ref{cor-main} (for typically real functions)}
Let $f=h+\overline{g} \in {\mathcal T}_H(K)$,  where $|\omega _{f}(z)|\leq k$ and $k=\frac{K-1}{K+1}$. In particular, $|\omega _{f}(0)|=|b_1|\leq k$. Clearly, each $f=h+\overline{g} \in {\mathcal T}_H(K)$ has the form
$$f=f_0+\overline{b_1f_0}=h_0+ \overline{b_1}g_0+\overline{g_0 + b_1h_0}
$$
where $f_0\in {\mathcal T}_H^0(K_0)$, $|\omega _{f_0}(z)|\leq k_0$ for all $z\in \ID$ and for some $k_0=\frac{K_0-1}{K_0+1}<1$, and $k_0=\frac{2k}{1+k^2}$.

Furthermore,  $f=f_0+\overline{b_1f_0}$ gives
$$ a_n(h)=a_n(h_0) +\overline{b_1}b_n(g_0) ~\mbox{ and }~b_n(g)=b_n(g_0) +b_1a_n(h_0).
$$
Finally, Theorem \ref{thm-main} and the fact that $|b_1|\leq k=\frac{1-\sqrt{1-k_0^2}}{k_0}$ give

\beqq
|a_n(h)| &\leq& A(n, k_0)+|b_1|B(n, k_0)\\
&\leq & A(n, k_0) + \frac{1-\sqrt{1-k_0^2}}{k_0}  B(n, k_0)=:E_n(k)\\
\eeqq
and
\beqq
|b_n(g)|&\leq& B(n, k_0)+|b_1|A(n, k_0)\\
&\leq & B(n, k_0) + \frac{1-\sqrt{1-k_0^2}}{k_0} A(n, k_0)=:F_n(k),
\eeqq
for $n=2, 3, \ldots$, where   $\phi_n(k_0)=A(n,k_0)$ and $\psi_n(k)=B(n,k_0)$ are as in Lemma \ref{lem1}  with $ k_0=\frac{2k}{1+k^2}$.
By Lemma \ref{lem1}, $A(n, k_0)$ and $B(n, k_0)$ are increasing functions of $k_0\in [0,1)$.  Also, it is easy to observe that $k_0\mapsto \frac{1-\sqrt{1-k_0^2}}{k_0}$ is increasing for $k_0\in (0,1)$,
and $k\mapsto \frac{2k}{1+k^2}$ is also increasing for $k\in [0,1)$. Consequently, both $E_n(k)$ and $F_n(k)$ are increasing functions of $k\in [0,1)$.



To see that the bounds are sharp, we only need consider functions $Q_k(z)$ given by \eqref{eq-qhkq},
since
$$\frac{2k}{1+k^2}=k_0~\Leftrightarrow~K^2=K_0.$$

The proof of the corollary for the other two cases follow similarly.
\qed
\medskip

 For the proofs of the other two cases of Theorem \ref{thm-main}, we need some useful results.

\begin{Lem}$($\cite[p.~108]{DuR}$)$\label{lemE}
If $f=h+\overline{g}\in {\mathcal S}_H$ is a starlike function, and if $H$ and $G$ are the analytic functions defined by
$$zH'(z)=h(z),\;\;zG'(z)=-g(z),\;\;H(0)=G(0)=0,$$
then $F=H+\overline{G}$ is a convex function of the class ${\mathcal C}_H$.
\end{Lem}

\begin{Lem}$($\cite{WangLZhang-01}$)$\label{lemF}
If $f=h+\overline{g}\in \mathcal{K}_H$, then there exist real numbers $\mu$, $\theta_0$ and a function $H(z)$ satisfying ${\rm Re}\,\{H(z)\}>0$ such that
$${\rm Re}\, \left\{H(z)\left[i e^{i \theta_0}(1-z^2)\big(e^{-i \mu} h^{\prime}(e^{i \theta_0} z)+e^{i \mu} g^{\prime}(e^{i \theta_0} z)\big)\right]\right\}>0
~\mbox{ for all $z\in \ID$}.$$
\end{Lem}

\subsection{Proof of  Theorem \ref{thm-main}}
{\bf Part 1:} First, we consider the case where ${\mathcal F}_0$ equals ${\mathcal S}_H^{*0}(K)$.
Let $f\in {\mathcal S}_H^{*0}$. Then by Lemma~F, 
the associated harmonic function $F=H+\overline{G}$ belongs to ${\mathcal C}_H^0$.
Appealing now to Lemma~D  
we conclude that there exist angles $\alpha$ and $\beta$ such that
$${\rm Re}\,\{(e^{i\alpha}H'(z)+e^{-i\alpha}G'(z))(e^{i\beta}-e^{-i\beta}z^2)\}>0 ~\mbox{ for all $z\in\ID$.}~
$$
Using the power series representation of $h$ and $g$ given by \eqref{111}, the last relation is equivalent to
${\rm Re}\,\{I(z)\}>0$, where
$$I(z)=\frac{e^{i\beta}-e^{-i\beta}z^2}{z}\sum_{n=1}^{\infty}(e^{i\alpha}a_n-e^{-i\alpha}b_n)z^n.
$$
Note that $a_1=1$ and $b_1=0$, since $f\in {\mathcal S}_H^0$. This last result says that
\be\label{eq-extra2}
e^{i\alpha}h(z)-e^{-i\alpha}g(z)=\sum_{n=1}^{\infty}(e^{i\alpha}a_n-e^{-i\alpha}b_n)z^n=\frac{z}{e^{i\beta}-e^{-i\beta}z^2}I(z).
\ee
As in the proof of Theorem \ref{thm1}, we see that the coefficients of the power series in the middle expression are dominated in modulus by those of
$$\frac{z}{1-z^2}\frac{1+z}{1-z}=\frac{z}{(1-z)^2}=\sum_{n=1}^{\infty}nz^n
$$
and this gives, $|e^{i\alpha}a_n-e^{-i\alpha}b_n|\leq n$  for $n\geq 1.$
Thus, $\big||a_n|-|b_n|\big|\leq n$ for $n\geq 1$.

For the other two estimates on the coefficients, we note that $g'(z)=\omega(z)h'(z)$, where $|\omega(z)|\leq k|z|$ in $\ID$.
Differentiation of \eqref{eq-extra2} gives
$$e^{i\alpha}h'(z)-e^{-i\alpha}g'(z)=\frac{d}{dz}\left\{\frac{z}{e^{i\beta}-e^{-i\beta}z^2}I(z)\right\},
$$
so that
$$g'(z)=\frac{\omega(z)}{e^{i\alpha}-e^{-i\alpha}\omega(z)}\frac{d}{dz}\left\{\frac{z}{e^{i\beta}-e^{-i\beta}z^2}I(z)\right\}.$$
This shows that the coefficients of $g'(z)$ are dominated in modulus by those of the function
$$\frac{kz}{1-kz}\Big(\frac{z}{(1-z)^2}\Big)'=\sum_{n =0}^{\infty}\left (\sum_{m =0}^{n}k^{m+1}(n-m+1)^2\right ) z^{n+1}.
$$
Finally, as in the proof of Theorem \ref{thm2}, one has the desired inequalities in this case.\\

{\bf Part 2:} Now, we consider the case where ${\mathcal F}_0$ equals ${\mathcal K}_H^{0}(K)$.
It follows from Lemmas E and G 
that the Taylor coefficients of the function
\beq\label{eq-F1}
F_1(z)=H(z)\left[i e^{i \theta_0}(1-z^2)\big(e^{-i \mu} h^{\prime}(e^{i \theta_0} z)+e^{i \mu} g^{\prime}(e^{i \theta_0} z)\big)\right]
\eeq
are dominated in modulus by the corresponding coefficients of the function $\frac{1+z}{1-z}$.
Let $F_2(z)=1/H(z)$. Then ${\rm Re}\,\{F_2(z)\}>0$ in $\ID$.
Hence, the coefficients of $e^{-i \mu} h^{\prime}(e^{i \theta_0} z)+e^{i \mu} g^{\prime}(e^{i \theta_0} z)$ are dominated in modulus by those of
$$\left(\frac{1+z}{1-z}\right)^2\frac{1}{1-z^2}=\frac{1+z}{(1-z)^3}.$$
Integration shows that the coefficients of $\int_0^{z}(e^{-i \mu} h^{\prime}(e^{i \theta_0} t)+e^{i \mu} g^{\prime}(e^{i \theta_0} t))\,dt$ are dominated
in modulus by those of the Koebe function $k_e(z)=z/(1-z)^2$.
Thus,
$$\big||a_n|-|b_n|\big|\leq |e^{-i\mu}a_n+e^{i\mu}b_n|\leq n ~\mbox{ for $n=2,3,\ldots$,}
$$
as the theorem asserts.

By \eqref{eq-F1}, we get
$$g^{\prime}(e^{i \theta_0} z)=-\frac{i e^{i(\mu-\theta_0)} F_1(z) F_2(z) \omega(e^{i \theta_0} z)}{(1-z^2)(1+e^{2 i \mu}\omega(e^{i \theta_0} z))}.$$
Again, the coefficients of $g'(e^{i \theta_0} z)$  are dominated in modulus by those of the function
$$\frac{kz}{1-kz}\frac{1+z}{(1-z)^3}=\sum_{n =0}^{\infty}\left (\sum_{m =0}^{n}k^{m+1}(n-m+1)^2\right) z^{n+1}.$$
Thus, as before, the inequalities
$$|b_{n}|\leq B(n,k)\;\;\mbox{and}\;\;|a_{n}|\leq A(n,k),
$$
hold for all $n\geq 2$. This completes the proof.
\qed

\subsection*{Acknowledgments}
The research was partly supported by the Natural Science Foundation of  China (No. 12371071).

\subsection*{Conflict of Interest Statement}
The authors declare that they have no conflict of interest, regarding the publication of this paper.

\subsection*{Data Availability Statement}
The authors declare that this research is purely theoretical and does not associate with any datas.



\begin{thebibliography}{99}

\bibitem{Ah} L. V. Ahlfors,
Lectures on quasiconformal mappings, Princeton (NJ), Van Nostrand company,
1966.

\bibitem{AAP-2019} R. M. Ali, Y. Abu Muhanna, and S. Ponnusamy,
{The spherical metric and univalent harmonic mappings},
\textit{Monatsh. Math.}  \textbf{188}(2019), 703--716.

\bibitem{Bo}
M. Borovikov, On Koebe radius and coefficients estimate for univalent harmonic
mappings, arXiv:2403.05526, (2024).

\bibitem{de_Branges} L.~de~Branges,
{A proof of the Bieberbach conjecture},
\textit{Acta Math.} {\bf 154}(1-2)(1985), 137--152.

\bibitem{BHH} D. Bshouty, W. Hengartner and O. Hossian,
Harmonic typically real mappings,
\textit{Math. Proc. Cambridge Philos. Soc.}, {\bf 119} (1996), 673--680.

\bibitem{CP-22} S. Chen and S. Ponnusamy, Koebe type theorems and pre-Schwarzian of harmonic
$K$-quasiconformal mappings, and their applications,
\textit{Acta Math. Sin. (Engl. Ser.)} \textbf{38}(2022), 1965--1980.

\bibitem{CSS} J. Clunie and T. Sheil-Small,
Harmonic univalent functions,
\textit{Ann. Acad. Sci. Fenn. Ser. A.I}, {\bf 9} (1984), 3--25.

\bibitem{Du} P. L. Duren,
Univalent functions,
Springer-Verlag, New York, 1983.

\bibitem{DuR} P. Duren,
Harmonic mappings in the plane,
Cambridge University Press, Cambridge, 2004.


\bibitem{LP-25} L. Li and S. Ponnusamy,
 Geometric subfamily of functions convex in some direction and
Blaschke products,
\textit{Bull. Malays. Math. Sci. Soc.} (2025), 12 pages;
 {\tt  DOI 10.1007/s40840-025-01867-9}

\bibitem{PonRasi2013} S. Ponnusamy and A. Rasila,
Planar harmonic and quasiregular mappings,
Topics in Modern  Function Theory: Chapter in
CMFT, RMS-Lecture Notes Series No. \textbf{19}(2013), 267--333.

\bibitem{SS} T. Sheil-Small,
Constants for planar harmonic mappings,
\textit{J. London. Math. Soc.}, {\bf 42} (1990), 237--248.

\bibitem{WangLZhang-01}
X. Wang, X. Liang and Y. Zhang,
{Precise coefficient estimates for close-to-convex harmonic univalent mappings},
\textit{J. Math. Anal. Appl.} {\bf 263}(2)(2001), 501--509

\bibitem{WLZ} X. Wang, X. Liang and Y. Zhang,
On harmonic typically real mappings,
\textit{J. Math. Anal. Appl.}, {\bf 277} (2003), 533--554.

\bibitem{wwrq} Z-G. Wang, X. Wang, A. Rasila and J. Qiu,
On a problem of Pavlovi$\acute{c}$ involving harmonic quasiconformal mappings,
https://doi.org/10.48550/arXiv.2405.19852.


\end{thebibliography}
\end{document}